\documentclass[12pt]{article}
\usepackage{amsmath,amsfonts,amssymb,amscd}
\usepackage{graphicx}

\newtheorem{theo}{Theorem}
\newtheorem{lem}{Lemma}[section]
\newtheorem{prop}{Proposition}[section]

\newtheorem{cor}{Corollary}[section]

\newtheorem{dfn}{Definition}[section]

\makeatletter \@addtoreset{equation}{section} \makeatother

\newcommand{\mC}{\mathbb{C}}

\newcommand{\mR}{\mathbb{R}}
\newcommand{\mT}{\mathbb{T}}
\newcommand{\mZ}{\mathbb{Z}}
\newcommand{\mN}{\mathbb{N}}

\newcommand{\bV}{{\bf V}}

\newcommand{\ba}{{\bf a}}
\newcommand{\bb}{{\bf b}}
\newcommand{\bc}{{\bf c}}

\newcommand{\calC}{{\cal C}}

\newcommand{\calE}{{\cal E}}

\newcommand{\calH}{{\cal H}}

\newcommand{\calK}{{\cal K}}

\newcommand{\calM}{{\cal M}}
\newcommand{\calN}{{\cal N}}

\newcommand{\calU}{{\cal U}}

\newcommand{\eps}{\varepsilon}
\newcommand{\ph}{\varphi}
\newcommand{\thet}{\vartheta}

\newcommand{\diag}{\operatorname{diag}}

\newcommand\qed{{\unskip\nobreak\hfil\penalty50
  \hskip2em\hbox{}\nobreak\hfil\mbox{\rule{1ex}{1ex} \qquad}
    \parfillskip=0pt \finalhyphendemerits=0\par\medskip}}

\begin{document}

%\large
\title
{On quantum Floquet theorem}
\author{D.~Treschev \\
Steklov Mathematical Institute of Russian Academy of Sciences
%\footnote{
%The research is supported by the RNF grant 14-50-00005.}
 }
\date{}
\maketitle

\begin{abstract}
We consider the Schr\"odinger equation $ih\partial_t\psi = H\psi$, $\psi=\psi(\cdot,t)\in L^2(\mT)$. The operator $H = -\partial^2_x + V(x,t)$ includes smooth potential $V$, which is assumed to be time $T$-periodic. Let $W=W(t)$ be the fundamental solution of this linear ODE system on $L^2(\mT)$. Then according to terminology from Lyapunov-Floquet theory, $\calM=W(T)$ is the monodromy operator. We prove that $\calM$ is unitarily conjugated to
$\exp\big(-\frac{T}{ih} \partial^2_x\big) + \calC$, where $\calC$ is a compact operator with an arbitrarily small norm.
\end{abstract}

\section{Introduction}

Let $V=V(x,t)$ be a smooth\footnote
{The corresponding functional spaces will be specified later.}
real-valued function, where $x\in\mT^d$, $\mT=\mR/ 2\pi\mZ$ and $t\in\mR$. Consider the Schr\"odinger equation
\begin{equation}
\label{Schrod}
  ih\partial_t\psi = H\psi, \qquad
  H = -\Delta + \eps V(x,t), \quad
  \Delta = \sum_{k=1}^d \partial^2_{x_k}.
\end{equation}
Here $\psi = \psi(x,t)$ is the wave function while $h>0$ is a constant. For any $t\in\mR$ the Schr\"odinger operator $H$ is an (unbounded) self-adjoint operator on $L^2(\mT^d)$. Hence the fundamental solution $W(t)$ (the Schr\"odinger propagator)
$$
  \psi(\cdot,0) \mapsto \psi(\cdot,t) = W(t)\psi(\cdot,0).
$$
is unitary for any real $t$.

The potential $V$ is assumed to be quasi-periodic in $t$. This means that there exists a smooth function $\bV:\mT^d\times\mT^n\to\mR$ and a constant vector $\omega\in\mR^n$ such that
$V(x,t)\equiv\bV(x,\omega t)$.

We are interested in the Floquet reducibility (FR) problem i.e., in the existence of
\begin{itemize}
\item a smooth quasiperiodic in time function $t\mapsto S(t)$ with values in the group of unitary operators on $L^2(\mT^d)$,
\item $t$-independent self-adjoint operator $\calH$ on $L^2(\mT^d)$
such that equation (\ref{Schrod}) is equivalent to
\begin{equation}
\label{calH}
  ih\partial_t\ph = \calH\ph, \qquad
  \psi = S\ph.
\end{equation}
\end{itemize}

It is known, \cite{Bam-Gra,Eli-Kuk}, that, if $\eps$ is sufficiently small and $V$ is real-analytic, then for most $\omega$ by using a quasi-periodic in time linear transformation $S$ of the Hilbert space $L^2(\mT^d)$ it is possible to reduce the equation (\ref{Schrod}) to the form (\ref{calH}), where $\calH = -\Delta + \eps Q$ and $Q$ is time independent bounded self-adjoint operator.\footnote
{The case $d=1$ is considered in \cite{Bam-Gra}, and the case $d\ge 2$ in \cite{Eli-Kuk}.}
The proof is based on methods of KAM-theory.

Note that this class of problems is close to the problems of reducibility of quasi-periodic cocycles. The corresponding references are contained in the survey paper \cite{Eli}.
\smallskip

In this paper we consider the case $d=n=\eps=1$. Hence we deal with the simplest version of the above FR problem except one aspect. We do not assume the potential $\eps V = V$ to be small.

The classical Lyapunov-Floquet approach to the FR problem for linear ODEs in a finite-dimensional space
\begin{equation}
\label{dotxi}
  \dot\xi = H(t) \xi, \qquad
  H(\cdot)\in L_m(\mC), \quad
  \xi\in\mC^m
\end{equation}
with $T$-periodic coefficients is based on construction of the monodromy operator $\calM = W(T)$, where $W(t)$ is the fundamental solution of (\ref{dotxi}).
Then (\ref{dotxi}) is equivalent to $\dot\eta=\calH\eta$, where $\calH$ may be computed from the equation $\calH = \frac1T \ln\calM$ and the linear $T$-periodic transformation $\xi=S(t)\eta$, is determined by $S(t) = W(t) e^{-t\calH}$.

To present our results, we start with some remarks.
\smallskip

{\bf 1}. Let $W = W(t)$ be the fundamental solution of equation (\ref{Schrod})$|_{d=n=\eps=1}$.
Since $W(t)$ is unitary for each\footnote
{Existence of solutions for $(\ref{Schrod})|_{d=n=\eps=1}$ will be proved below.}
$t\in\mR$, the monodromy operator $\calM=W(T)$ is also unitary.
\smallskip

{\bf 2}. In the case $n=1$ w.l.g. $\omega=\omega_1=1$. The general case is reduced to this one by a change of time. Hence, below $V$ is assumed to be $2\pi$-periodic in $t$ and $T=2\pi$.
\smallskip

{\bf 3}. Let $v_0(t) = \frac1{2\pi}\int_\mT V(x,t)\, dx$ be the space average of the potential and
$v_{00} = \frac1{2\pi}\int_\mT v_0(t)\, dt$. By using $2\pi$-periodic in $t$ change of the unknown function $\psi\mapsto\tilde\psi$,
$$
  \psi(x,t) = \tilde\psi(x,t) \exp\Big(\frac1{ih}\int_0^t (v_0(\tau) - v_{00})\, d\tau\Big),
$$
it is possible to replace in $V$ the space average $v_0(t)$ by $v_{00}$.

Moreover, the operators $W(t)$, corresponding to the potential $V$ and $V - v_{00}$ differ by the scalar multiplier $e^{\frac{v_{00}}{ih}t}$. Hence, it is sufficient to study the case $v_{00}=0$. Below we assume $v_0(t)\equiv 0$.
\smallskip

{\bf 4}. Let $U$ be a (constant) invertible operator on $L^2(\mT)$. If $\psi=U\tilde\psi$ is a coordinate change, then in the new coordinates the fundamental solution takes the form $\tilde W = U^{-1}WU$. In particular, the new monodromy operator is $\tilde\calM = U^{-1}\calM U$. Due to this $\calM$ is assumed to be defined up to a conjugacy.
\smallskip

Consider the standard Fourier basis $e^{ikx}$, $k\in\mZ$ in $L^2(\mT)$. In this basis $H=A+B$,
\begin{equation}
\label{coord}
  A = \diag(k^2), \quad
  B = (B_{jk}), \qquad
  B_{jk} = v_{j-k},\quad
  V = \sum_{k\in\mZ} v_k(t) e^{ikx}.
\end{equation}
For any $t\in\mT$ the infinite vector $(v_k(t))_{k\in\mZ}$ is an element of the Hilbert space $l^2=l^2(\mZ)$, dual to $L^2(\mT)$.

Equation (\ref{Schrod}) is the following ODE system on $l^2(\mZ)$:
\begin{equation}
\label{Schrod1}
     ih\dot\psi_k
  =  k^2\psi_k + \sum_{j\in\mZ} v_{k-j}(t) \psi_j, \qquad
     \psi = \sum_{k\in\mZ} \psi_k e^{ikx}.
\end{equation}

Now we specify the smoothness classes for the potential.

\begin{dfn}
We say that $V:\mT^2\to\mR$ lies in $C^{\alpha,\beta,\gamma}$, $\alpha,\gamma\in\mN$, $\beta\ge 0$ if $V=V(x,t)$ is $C^\gamma$-smooth in $t$ and its Fourier expansion (\ref{coord}) satisfies the estimate
\begin{equation}
\label{Vsmooth}
  \|v_k\|_{C^\gamma} \le \frac{c_v e^{-\beta|k|}}{\langle k\rangle^\alpha}, \qquad
  \langle k\rangle = |k| + 1.
\end{equation}
\end{dfn}

The case $\beta=0$ is close to $C^\alpha$-smoothness in $x$. If $\beta>0$ then $V$ is real-analytic in $x$.

\begin{theo}
\label{theo:mono}
Suppose $V\in C^{\alpha,\beta,2}$, where $\alpha>2$. Then the operator $\calM$ is unitarily conjugated to $\calM^{(0)} + \calM^{(d)} + \calM^{(c)}$, where
\begin{eqnarray}
\label{M0}
     \calM^{(0)}
 &=& \exp\Big(-\frac{2\pi}{ih}\partial^2_x\Big), \quad
     \calM^{(0)}_{jk}
  =  \delta_{jk} e^{\frac{2\pi}{ih} k^2}, \\
\label{Md}
     \calM^{(d)}_{jk}
 &=& \delta_{-jk} \frac{e^{\frac{2\pi}{ih}k^2}}{ih} \int_0^{2\pi} v_{-2k}(\tau)\, d\tau , \\
\label{Mc}
     |\calM^{(c)}_{jk}|
 &\le& \frac{c e^{-\beta |j-k|}}
            {\langle j\rangle \langle k\rangle \langle j-k\rangle^{\alpha-1}}.
\end{eqnarray}
The constant $c$ depends only on $h,c_v$, and $\alpha$.
\end{theo}

Estimates (\ref{Mc}) mean that $\calM^{(c)}$ is a compact operator. The estimate (\ref{Vsmooth}) admits the possibility
$$
  |\calM^{(d)}_{-kk}| \ge \frac{\hat c_v e^{-2\beta |k|}}{\langle k\rangle^\alpha}, \qquad
  \hat c_v > 0.
$$
This is much greater than $|\calM^{(c)}_{-kk}|$ for large $|k|$. Hence we may regard $\calM^{(d)}$ as a ``more important correction'' to $\calM^{(0)}$ in comparison with $\calM^{(c)}$.

It is possible to make the monodromy operator $\calM^{(0)} + \calM^{(d)} + \calM^{(c)}$ arbitrarily close to a diagonal one by using another unitary conjugacy. Here is the result.

\begin{theo}
\label{theo:mono1}
Suppose $V$ is $C^{\alpha,\beta,2}$, $\alpha>2$. Then for any $N\ge N_0(h,c_v,\alpha)$ the operator $\calM$ is unitarily conjugated to $\tilde\calM^{(0)} + \tilde\calM^{(d)} + \tilde\calM^{(c)}$, where\footnote
{The quantities $d_k$, $|k|\le N$ will be determined in Section \ref{sec:M+M}. }
\begin{eqnarray}
\label{tildeM}
     \tilde\calM^{(0)}_{jk}
 &=& d_j\delta_{jk}, \qquad
     |d_j|=1, \quad d_k = e^{\frac{2\pi}{ih}k^2}\quad
     \mbox{for any } |k| > N,  \\
\label{tildeMd}
     \tilde\calM^{(d)}_{jk}
 &=& \left\{\begin{array}{cl}
                     0       &  \mbox{if}\quad  -N\le j,k\le N, \\
            \calM^{(d)}_{jk} &  \mbox{if \quad  $|j|>N$ or $|k|>N$}, 
            \end{array}
     \right. \\
\label{tildeMc}
       |\tilde\calM^{(c)}_{jk}|
 &\le& \left\{\begin{array}{cl}
       \tilde c (N+1)^{-2} & \mbox{ if }\; -N\le j,k\le N, \\[1mm]
       \tilde c e^{-2\beta (|j|-N)}\langle j\rangle^{-2}
                           & \mbox{ if }\quad |k|\le N < |j| , \\[1mm]
       \tilde c e^{-2\beta (|k|-N)}\langle k\rangle^{-2}
                           & \mbox{ if }\quad |j|\le N < |k| , \\[1mm]
\displaystyle
       \frac{\tilde c e^{-2\beta (|j-k|)}}
            {\langle j\rangle \langle k\rangle \langle j-k\rangle^{\alpha-1}}
                           & \mbox{ if }\quad |j|> N$ and $|k|>N ,
             \end{array}
      \right.
\end{eqnarray}
$\tilde c = \tilde c(h,c_v,\alpha)$.
\end{theo}

Hence $\calM$ is unitarily conjugated to an operator which differs from a diagonal one by a compact operator with an arbitrarily small norm. Computation of the monodromy operator is the first step in the Lyapunov-Floquet reducibility procedure.
However computation of the operator
$\calH = \frac1{2\pi} \ln (\tilde\calM^{(0)}+\tilde\calM^{(d)}+\tilde\calM^{(c)})$
remains a nontrivial problem which requires application of methods of KAM-theory. Iteration procedures for construction of $\calH$ include small divisors close to
$e^{\frac{2\pi}{ih} k_1^2} - e^{\frac{2\pi}{ih} k_2^2}$, $\, k_1,k_2\in\mZ$.

Plan of the paper is as follows. In Section \ref{sec:mono} we present the monodromy operator as the sum of four infinite matrices, where the first three, $\calM^{(0)},\calM^{(d)}$, and $\calM^{(1)}$, are explicit and the last one, $\calM^{(2)}$, is ``small'' (Theorem \ref{theo:approx}). Then by using a unitary conjugacy, we remove $\calM^{(1)}$. This gives a proof of Theorem \ref{theo:mono}.

In Section \ref{sec:est} we prove Theorem \ref{theo:approx} by studying the fundamental solution $W(t)$ in the Fourier basis. Here we use the fact that $k_0$-th column of the infinite matrix $W(t)$ is the solution $\psi$ of equation (\ref{Schrod1}) with the initial condition $\psi_k(0)=\delta_{kk_0}$.

Section \ref{sec:diag} contains proof of Theorem \ref{theo:mono1}: we take the operator $\calM^{(0)}+\calM^{(d)}+\calM^{(c)}$ from Theorem \ref{theo:mono} and diagonalize (approximately) its middle $(2N+1)\times (2N+1)$ block by a unitary conjugacy. Section \ref{sec:tech} contains auxiliary technical statements and a somewhat cumbersome proof of Lemma \ref{lem:Gn}.

\section{Monodromy operator}
\label{sec:mono}

Let $\calM = W(2\pi)$ be the monodromy operator. Since $\calM=(\calM^*)^{-1}$ is unitary, we have:
$\calM_{jk} = \big(\overline{\calM^{-1}}\big)_{kj}$.

Consider infinite matrices $\calM^{(0)}$, $\calM^{(d)}$ (see (\ref{M0})--(\ref{Md})), and $\calM^{(1)}$, where
$$
    \calM^{(1)}_{jk}
  = \left\{\begin{array}{cc}
           0                &     \mbox{if } j^2 = k^2, \\
    \displaystyle
    \frac{v_{j-k}(0) ( e^{\frac{T}{ih} k^2} - e^{\frac{T}{ih} j^2} ) }
         {k^2-j^2}          &     \mbox{if } j^2\ne k^2.
           \end{array}
    \right.
$$

\begin{theo}
\label{theo:approx}
There exists a constant $c_\calM$ depending only on $c_v,h$, and $\alpha$ such that for all $j,k\in\mZ$
\begin{equation}
\label{calM}
         \calM
     =   \calM^{(0)} + \calM^{(d)} + \calM^{(1)} + \calM^{(2)}, \quad
         |\calM^{(2)}_{jk}|
    \le  \frac{c_\calM e^{-\beta |j-k|} }
              {\langle j\rangle \langle k\rangle \langle j-k\rangle^{\alpha-1}}.
\end{equation}
\end{theo}

We prove Theorem \ref{theo:approx} in Section \ref{sec:est}.
\smallskip

Consider the anti-selfadjoint operator
$$
  G = (G_{jk}), \quad
    G_{jk}
  = \left\{ \begin{array}{cc}\displaystyle
                                  0                      & \mbox{if }  j^2 = k^2, \\
                             \displaystyle
                              \frac{v_{j-k}(0)}{j^2-k^2} & \mbox{if }  j^2\ne k^2.
            \end{array}
    \right.
$$

\begin{lem}
\label{lem:Gn}
For any $n\ge 2$
\begin{equation}
\label{Gn<}
      |G_{jk}|
  \le \frac{3c_v e^{-\beta |j-k|}}{(\langle j\rangle+\langle k\rangle) \langle j-k\rangle^\alpha}, \quad
      |(G^n)_{jk}|
  \le \frac{(3c_v)^n c_\alpha^{n-1} e^{-\beta |j-k|}}
           {\langle j\rangle \langle k\rangle \langle j-k\rangle^\alpha},
\end{equation}
where $c_\alpha=c_\nu|_{\nu=\alpha}$ is the constant from Lemma \ref{lem:tech1}.
\end{lem}

\begin{cor}
\begin{equation}
\label{calM1}
      |\calM^{(1)}_{jk}|
  \le \frac{6 c_v e^{-\beta |j-k|}}
           {(\langle j\rangle + \langle k\rangle) \langle j-k\rangle^\alpha}.
\end{equation}
\end{cor}

We prove Lemma \ref{lem:Gn} in Section \ref{sec:Gn}. Lemma \ref{lem:Gn} implies that
\begin{equation}
\label{e^G}
  e^{\pm G} = I \pm G + G^\pm, \qquad
  |G^\pm_{jk}| \le \frac{c_e e^{-\beta |j-k|}}
                        {\langle j\rangle \langle k\rangle \langle j-k\rangle^\alpha}, \quad
  c_e = c_e(c_v,\alpha).
\end{equation}

Now Theorem \ref{theo:mono} is a direct consequence from the following

\begin{cor}
\label{cor:N}
The unitary operator $S = e^{-G}$ conjugates $\calM$ with $\calN$:
\begin{equation}
\label{SMS}
  \calN = S^* \calM S, \quad
      \big| \calN_{jk} - \calM^{(0)}_{jk} -  \calM^{(d)}_{jk} \big|
  \le \frac{c e^{-\beta |j-k|}}
           {\langle j\rangle \langle k\rangle \langle j-k\rangle^{\alpha-1}}, \qquad
      c = c(c_v,h,\alpha).
\end{equation}
\end{cor}

To prove (\ref{SMS}), we note that
$e^G\calM e^{-G} = \calM^{(0)} + \calM^{(d)} + \calK^{(1)} + \calK^{(2)}$,
\begin{eqnarray*}
     \calK^{(1)}
 &=& \calM^{(1)} - \calM^{(0)} G + G \calM^{(0)} \, =\, 0, \\
     \calK^{(2)}
 &=& \calM^{(2)} + \calM^{(0)} G^- + G^+ \calM^{(0)} \\
 &&  + (\calM^{(d)} + \calM^{(1)} + \calM^{(2)}) (-G + G^-)
     +\, (G + G^+) (\calM^{(d)} + \calM^{(1)} + \calM^{(2)}) \\
 &&  + (G + G^+) (\calM^{(0)} + \calM^{(d)} + \calM^{(1)} + \calM^{(2)}) (-G + G^-).
\end{eqnarray*}

The operator $\calM^{(2)}$ in $\calK^{(2)}$ is estimated in (\ref{calM}), $\calM^{(0)} G^-$ and $G^+\calM^{(0)}$ satisfy the same estimates as $G^\pm$, see (\ref{e^G}). By (\ref{calM}) and (\ref{calM1})
$$
      \big| \calM^{(d)}_{jk} + \calM_{jk}^{(1)} + \calM^{(2)}_{jk} \big|
  \le \frac{\hat c e^{-\beta |j-k|}}
           {(\langle j\rangle + \langle k\rangle) \langle j-k\rangle^{\alpha-1}}, \qquad
      \hat c = \hat c(c_v,\alpha,h).
$$
This inequality combined with Lemma \ref{lem:tech1} imply estimates for remaining terms in $\calK^{(2)}$. \qed

\section{Estimates for solutions of (\ref{Schrod1})}
\label{sec:est}

Let $\psi$ be the $k_0$-th column of the operator $W$ in the Fourier basis. Then $\psi$ is the solution of the initial value problem
\begin{equation}
\label{kk0}
  ih\dot\psi_k = k^2\psi_k + \sum_{j\in\mZ} v_{k-j}(t) \psi_j, \qquad
  \psi_k(0) = \delta_{kk_0}.
\end{equation}
We assume $|k_0|$ to be big. This will be our large parameter. As a first approximation for the solution of (\ref{kk0}) we take $\psi^{(1)}$, the solution of the initial value problem
\begin{equation}
\label{kk0(1)}
  ih\dot\psi_k^{(1)} = k^2\psi_k^{(1)} + v_{k-k_0}(t) \psi_{k_0}^{(1)}, \qquad
  \psi_k(0) = \delta_{kk_0}.
\end{equation}
This problem is obtained from (\ref{kk0}) if one takes in the sum only the term, corresponding to $j=k_0$. Explicit solution of (\ref{kk0(1)}) (for $k\ne\pm k_0$) is
\begin{eqnarray}
\nonumber
     \psi_{k_0}^{(1)}(t)
 &=& e^{\frac{t}{ih} k_0^2}, \qquad
      \psi_{-k_0}^{(1)}(t)
\,=\, \frac{e^{\frac{t}{ih} k_0^2}}{ih} \int_0^t v_{-2k_0}(\tau)\, d\tau, \\
\nonumber
     \psi_k^{(1)}(t)
 &=& \frac{e^{\frac{t}{ih}k^2}}{ih} \int_0^t v_{k-k_0}(\tau) e^{\frac{\tau}{ih}(k_0^2-k^2)}\, d\tau
  =  \frac{v_{k-k_0}(t) e^{\frac{t}{ih} k_0^2} - v_{k-k_0}(0) e^{\frac{t}{ih} k^2} }
          {k_0^2-k^2}
    + \psi_k^{(11)}(t), \\
\label{psi1}
     \psi_{\pm k_0}^{(11)}(t)
 &=& 0, \qquad
        \psi_k^{(11)}(t)
 \, =\, - e^{\frac{t}{ih} k^2} 
            \int_0^t \frac{\dot v_{k-k_0}(\tau)}{k_0^2-k^2}
                     e^{\frac{\tau}{ih} (k_0^2-k^2)} \, d\tau .
\end{eqnarray}
Here we have used integration by parts.

\begin{prop}
\label{prop:approx}
Suppose $|t|\le t_0$. Then for some constant $c_\psi$ depending only on $c_v,h,\alpha$, and $t_0$
\begin{equation}
\label{psi-psi}
       |\psi_k(t) - \psi_k^{(1)}(t)|
   \le \frac{c_\psi e^{-\beta |k-k_0|}}
            {\langle k_0\rangle^2 \langle k-k_0\rangle^{\alpha-1}} .
\end{equation}
\end{prop}

{\it Proof of Proposition \ref{prop:approx}}. We put
$\thet_k = e^{-\frac{t}{ih} k^2} (\psi_k - \psi_k^{(1)})$. Then by (\ref{kk0}) and (\ref{kk0(1)})
\begin{eqnarray}
\nonumber
     ih\dot\thet_k
 &=& \sum_{j\in\mZ} e^{\frac{t}{ih} (j^2-k^2)} v_{k-j} \thet_j + \dot\Phi_k,
     \qquad  \thet_k(0) = 0, \\
\label{Phi<}
     \Phi_k(t)
 &=& \int_0^t
        e^{-\frac{\tau}{ih} k^2}
          \sum_{k_0\ne j\in\mZ} v_{k-j}(\tau) \psi_j^{(1)}(\tau)  \, d\tau.
\end{eqnarray}

\begin{lem}
\label{lem:Phi}
Suppose $|t|\le t_0$. Then for some constant $c_\Phi = c_\Phi(h,c_v,\alpha,t_0)$
$$
      |\Phi_k(t)|
  \le \frac{c_\Phi e^{-\beta |k-k_0|}}
           {\langle k_0\rangle^2 \langle k-k_0\rangle^{\alpha-1}} .
$$
\end{lem}

Proof of Lemma \ref{lem:Phi} is straightforward, but cumbersome. We present it in Section \ref{sec:Phi}.
\smallskip

We denote
$$
    \xi_k
  = \Big(\thet_k - \frac1{ih}\Phi_k\Big) \langle k-k_0\rangle^{\alpha-1} e^{\beta |k-k_0|}.
$$
Then
\begin{eqnarray*}
     \frac{ih e^{-\beta |k-k_0|}}{\langle k-k_0\rangle^{\alpha-1}} \dot\xi_k
 &=& \sum_{j\in\mZ} e^{\frac{t}{ih}(j^2-k^2)}
                    \frac{v_{k-j} e^{-\beta |j-k_0|}}
                         {\langle j-k_0\rangle^{\alpha-1}} \xi_j
    + \Psi_k(t), \qquad
     \xi_k(0) = 0, \\
     \Psi_k(t)
 &=& \sum_{j\in\mZ} \frac{e^{\frac{t}{ih}(j^2-k^2)}}{ih}  
                    v_{k-j} \Phi_j(t) .
\end{eqnarray*}
By Lemmas \ref{lem:Phi} and \ref{lem:tech1} for $|t|\le t_0$
$$
      |\Psi_k(t)|
  \le \sum_{j\in\mZ} \frac{h^{-1} c_v c_\Phi e^{-\beta |k-k_0|}}
                          {\langle k-j\rangle^\alpha \langle k_0\rangle^2
                               \langle j-k_0\rangle^{\alpha-1}}
  \le \frac{h^{-1} c_v c_\Phi c_{\alpha-1} e^{-\beta |k-k_0|}}
           {\langle k_0\rangle^2 \langle k-k_0\rangle^{\alpha-1}}.
$$
We obtain:
$$
     h|\dot\xi_k|
 \le \sum_{j\in\mZ} \frac{c_v \langle k-k_0\rangle^{\alpha-1}}
                         {\langle k-j\rangle^\alpha \langle j-k_0\rangle^{\alpha-1}} |\xi_j|
     + \frac{c_v c_\Phi c_{\alpha-1}}{h \langle k_0\rangle^2}, \qquad
     \xi_k(0) = 0.
$$

By Lemma \ref{lem:tech1}
$$
  \sum_{j\in\mZ} \frac{\langle k-k_0\rangle^{\alpha-1}}
                      {\langle k-j\rangle^\alpha \langle j-k_0\rangle^{\alpha-1}}
  \le c_{\alpha-1}.
$$
This implies that for all $k\in\mZ$ we have: $|\xi_k(t)|\le |\mu(t)|$, where
$$
     h\dot\mu
  =  c_v c_{\alpha-1} \mu + \frac{c_v c_\Phi c_{\alpha-1}}{h \langle k_0\rangle^2}, \qquad
     \mu(0) = 0.
$$
Solving this equation, we obtain the estimate
$$
      |\xi_k(t)|
  \le \mu(t)
  \le \frac{c_\xi}{\langle k_0\rangle^2}, \qquad
  c_\xi = \frac{c_\Phi}{h} \Big(e^{\frac{c_v c_{\alpha-1}}{h} t} - 1\Big).
$$
Hence
$$
     |\psi_k - \psi_k^{(1)}|
  =  |\thet_k|
 \le \frac{|\Phi_k|}{h} + \frac{|\xi_k| e^{-\beta |k-k_0|}}
          {\langle k-k_0\rangle^{\alpha-1}}
 \le \frac{(h^{-1} c_\Phi + c_\xi) e^{-\beta |k-k_0|}}
          {\langle k_0\rangle^2 \langle k-k_0\rangle^{\alpha-1}} .
$$
Proposition \ref{prop:approx} is proved.  \qed

Now we prove Theorem \ref{theo:approx}. We take in Proposition \ref{prop:approx} $t_0=2\pi$. Then
$\psi_k(2\pi) = \calM_{kk_0}$. Therefore
$\psi_k(2\pi)
  = \calM^{(0)}_{kk_0} + \calM^{(d)}_{kk_0} + \calM^{(1)}_{kk_0} + \calM^{(2)}_{kk_0}$,
where $\calM^{(2)}_{kk_0} = \psi_{kk_0} - \psi_{kk_0}^{(1)} + \psi_k^{(11)}$. We estimate $\psi_k^{(11)}$ ($k\ne\pm k_0$) by using integration by parts:
$$
     |\psi_k^{(11)}|
  =  \bigg| \int_0^t \frac{h \dot v_{k-k_0}(\tau)}{(k_0^2-k^2)^2}\, de^{\frac{\tau}{ih}(k_0^2-k^2)}
     \bigg|
 \le \frac{|t| h c_v e^{-\beta |k-k_0|}}{\langle k_0\rangle^2 \langle k-k_0\rangle^\alpha}.
$$
this estimate combined with (\ref{psi-psi}) implies
\begin{equation}
\label{calM2}
      |\calM^{(2)}_{kk_0}|
  \le \frac{c_\calM e^{-\beta|k-k_0|}}{\langle k_0\rangle^2 \langle k-k_0\rangle^{\alpha-1}}.
\end{equation}

By the same argument we obtain:
$$
    W(-2\pi)
  = \calM^{-1}
  = \calM^{(0)-}_{kk_0} + \calM^{(d)-}_{kk_0} + \calM^{(1)-}_{kk_0} + \calM^{(2)-}_{kk_0},
$$
where $\calM^{(0)-}_{kk_0}, \calM^{(d)-}_{kk_0}, \calM^{(1)-}_{kk_0}, \calM^{(2)-}_{kk_0}$ are defined in the same way as $\calM^{(0)}_{kk_0}, \calM^{(d)}_{kk_0}, \calM^{(1)}_{kk_0}, \calM^{(2)}_{kk_0}$, but instead of $t=2\pi$ one should take $t=-2\pi$. Then
\begin{equation}
\label{calM2-}
      |\calM^{(2)-}_{kk_0}|
  \le \frac{c_\calM e^{-\beta|k-k_0|}}{\langle k_0\rangle^2 \langle k-k_0\rangle^{\alpha-1}}.
\end{equation}

By definition of $\calM^{(0)}_{kk_0}, \calM^{(d)}_{kk_0}$, and $\calM^{(1)}_{kk_0}$, we have:
$$
  \calM^{(0)-}_{kk_0} = \big(\calM^{(0)}_{kk_0}\big)^*, \quad
  \calM^{(d)-}_{kk_0} = \big(\calM^{(d)}_{kk_0}\big)^*, \quad
  \calM^{(1)-}_{kk_0} = \big(\calM^{(1)}_{kk_0}\big)^*.
$$
Therefore $\calM^{(2)-}_{kk_0} = \big(\calM^{(2)}_{kk_0}\big)^*$ and by (\ref{calM2-})
\begin{equation}
\label{calM2*}
      |\calM^{(2)}_{kk_0}|
   =  |\calM^{(2)-}_{k_0k}|
  \le \frac{c_\psi e^{-\beta|k-k_0|}}{\langle k\rangle^2 \langle k-k_0\rangle^{\alpha-1}}.
\end{equation}
Estimates (\ref{calM2}) and (\ref{calM2*}) imply (\ref{calM}). \qed

\section{Approximate diagonalization}
\label{sec:diag}

In this section we prove Theorem \ref{theo:mono1}. More precisely, we consider the operator $\calN$ from Corollary \ref{cor:N} and construct a unitary $\calU$ such that
$\calU^{-1}\calN\calU - \tilde\calM^{(0)}$ is a compact operator with an arbitrarily small norm.

\subsection{Step 1}

Take $\calN = \calM^{(0)} + \calM^{(d)} + \calM^{(c)}$ as in Theorem \ref{theo:mono}.
Let $w_{-N}, \ldots, w_N$ be columns of the $(2N+1)\times(2N+1)$ matrix $(\calN_{jk})_{-N\le j,k\le N}$. Then
\begin{equation}
\label{w-delta}
      \bigg| w_{jk} - \delta_{jk} e^{\frac{2\pi}{ih} k^2}
                    - \delta_{-jk} \frac{e^{\frac{2\pi}{ih}k^2}}{ih} \int_0^{2\pi} v_{-2k}(\tau)\, d\tau
      \bigg|
  \le \frac{c e^{-\beta |j-k|}}
           {\langle j\rangle \langle k\rangle \langle j-k\rangle^{\alpha-1}}.
\end{equation}

The operator $\calN$ is unitary. Hence the scalar products $\langle w_{j'},w_{j''}\rangle$ satisfy the estimates
$$
       \langle w_{j'},w_{j''}\rangle - \delta_{j' j''}
   =   \sum_{|k|\le N} \calN_{kj'} \overline \calN_{kj''} - \delta_{j'j''}
   =   \sum_{|k| > N} w_{kj'} \overline w_{kj''} .
$$
In this sum $|k|>N$ while $|j'|\le N$ and $|j''|\le N$. Therefore
$\delta_{\pm kj'}=\delta_{\pm k j''}=0$. Hence by (\ref{w-delta})
\begin{eqnarray}
\nonumber
       |\langle w_{j'},w_{j''}\rangle - \delta_{j' j''}|
 &\le& \sum_{|k|>N} \frac{c^2 e^{-\beta (|j'-k|+|j''-k|)}}
                        {\langle j'\rangle \langle j''\rangle \langle k\rangle^2
                          \langle j'-k\rangle^{\alpha-1} \langle j''-k\rangle^{\alpha-1}}
 \,\le\,   \eps_{j'j''},  \\
\label{eps<}
       \eps_{j'j''}
  &=&  \frac{c^2 c_{\alpha-1} e^{-\beta |j'-j''|}}
           {\langle j'\rangle \langle j''\rangle (N+1)^2 \langle j'-j''\rangle^{\alpha-1}}.
\end{eqnarray}

\begin{prop}
\label{prop:U}
For any sufficiently large $N\in\mN$ there exists $U=(u_{jk})\in U(2N+1)$ such that the operator
$$
  W' = (w'_{jk}) = (u_{jk}-w_{jk}), \qquad
  -N\le j,k\le N.
$$
satisfies the estimate
\begin{equation}
\label{U-N}
        |w'_{jk}|
   \le  \sigma_k |w_{jk}| + \sum_{-|k|\le l < |k|,\, l\ne k} 4\eps_{lk} |w_{jl}|, \qquad
        \sigma_k
   =    4\sum_{|l|\le |k|} \eps_{lk}.
\end{equation}
\end{prop}

\begin{cor}
For some constant $c_w=c_w(c_v,h,\alpha)$
\begin{equation}
\label{w'}
       |w'_{jk}|
  \le  \frac{c_w e^{-\beta |j-k|}}
            {\langle j\rangle \langle k\rangle (N+1)^2 \langle j-k\rangle^{\alpha-1}}, \qquad
       -N\le j,k\le N.
\end{equation}
\end{cor}

Indeed, by (\ref{eps<})
$$
      |\sigma_k|
  \le 4\sum_{|l|\le |k|}
           \frac{c^2_\calN c_{\alpha-1} e^{-\beta |l-k|}}
                {\langle l\rangle \langle k\rangle (N+1)^2 \langle l-k\rangle^{\alpha-1}}
  \le \frac{c^2_\calN c_{\alpha-1} c_*}
           {\langle k\rangle^2 (N+1)^2},
$$
where $c_*$ is an absolute constant. Hence
\begin{eqnarray*}
       \sigma_k |w_{jk}|
 &\le& \bigg( \delta_{jk}
            + \delta_{-jk} \frac{2\pi}{h}
              \frac{c_v e^{-\beta |k|}}
                   {\langle 2k\rangle^\alpha}
            + \frac{c e^{-\beta |j-k|}}
                   {\langle j\rangle \langle k\rangle \langle j-k\rangle^{\alpha-1}} \bigg)
       \frac{c^2 c_{\alpha-1} c_*}
            {\langle k\rangle^2 (N+1)^2} \\
 &\le& \frac{c_\sigma e^{-\beta |j-k|}}
            {\langle j\rangle \langle k\rangle \langle j-k\rangle^{\alpha-1} (N+1)^2}.
\end{eqnarray*}
We also have:
\begin{eqnarray*}
       \sum_{|l|\le |k|,\, l\ne k} 4\eps_{lk} |w_{jl}|
 &\le& \sum_{|l|\le |k|,\, l\ne k}
         \frac{4c^2 c_{\alpha-1} e^{-\beta |l-k|}}
              {\langle l\rangle \langle k\rangle \langle l-k\rangle^{\alpha-1} (N+1)^2}
      \times \\
 &&\qquad\qquad \times
         \bigg( \delta_{jl}
            + \delta_{-jl} \frac{2\pi}{h}
              \frac{c_v e^{-\beta |l|}}
                   {\langle 2l\rangle^\alpha}
            + \frac{c e^{-\beta |j-l|}}
                   {\langle j\rangle \langle l\rangle \langle j-l\rangle^{\alpha-1}} \bigg)   \\
 &\le& \frac{4c^3 c_{\alpha-1} c_\eps e^{-\beta |j-k|}}
            {\langle j\rangle \langle k\rangle \langle j-k\rangle^{\alpha-1} (N+1)^2}, \qquad
       c_\eps = c_\eps(h,c_v,\alpha).     
\end{eqnarray*}
This implies (\ref{w'}). \qed

We start proof of Proposition \ref{prop:U} from the following simple lemma.

\begin{lem}
\label{lem:eee}
If $N$ is sufficiently large,
$$
  \sum_{l=-N}^N \eps_{jl}\eps_{lk} \le \frac12 \eps_{jk} \quad
  \mbox{for any }  j,k\in\{-N,\ldots,N\}.
$$
\end{lem}

{\it Proof}. Direct computation with the help of the inequality (\ref{eps<}). \qed

We have to perform an orthogonalization of the ``almost Hermitian'' basis $w_{-N},\ldots,w_N$. The standard Gram-Schmidt orthogonalization procedure gives bad estimates for the differences $U_{jk}-w_{jk}$. We use another version of the orthogonalization process, based on the following lemma.

\begin{lem}
\label{lem:GS}
Let $a_1,\ldots,a_s,a_{s+1}$ be linearly independent vectors such that
\begin{equation}
\label{aEe}
  \langle a_j,a_k\rangle = \delta_{jk} + \calE_{jk}, \quad
  |\calE_{jk}| \le \sigma_{jk}, \quad
  \sum_{l=1}^s \sigma_{jl} \sigma_{lk} \le \frac12 \sigma_{jk} \quad
  \mbox{for any } j,k\in \{1,\ldots,s\}.
\end{equation}
Then there exists a vector
$\displaystyle b=(1+\sigma')^{-1} \big( a_{s+1} - \sum_{j=1}^s \lambda_j a_j \big)$ such that
\begin{eqnarray}
\nonumber
&    |b|=1, \quad
     \langle a_j,b\rangle = 0, \qquad  1\le j\le s,  & \\
\label{b-a}
&\displaystyle
       |\lambda_j| \le 2\sigma_{j s+1}, \quad
       |\sigma'|
  \le  \sigma_{s+1 s+1} + \sum_{j=1}^s  2\sigma_{j s+1} (1 + \sigma_{jj}). &
\end{eqnarray}
\end{lem}

{\it Proof}. First note that the equation $(|a_j|-1)(|a_j|+1) = \sigma_{jj}$ implies
\begin{equation}
\label{|a-1|}
  \big| |a_j| - 1 \big| \le \sigma_{jj}, \qquad  j=1,\ldots,s+1.
\end{equation}

We put $b' = a_{s+1} - \sum_{j=1}^s \lambda_j a_j$ and assume that
$\langle a_k,b'\rangle=0$ for any $k=1,\ldots,s$. Then
$$
  \langle a_k,a_{s+1} \rangle = \sum_{j=1}^s \langle a_k,a_j\rangle \lambda_j.
$$
This implies $(I+\calE)\lambda = \mu$, where
$$
  \lambda = \left(\begin{array}{c} \lambda_1 \\ \ldots \\ \lambda_s
                  \end{array}\right), \quad
  \mu = \left(\begin{array}{c} \langle a_1,a_{s+1} \rangle \\
                                  \ldots \\
                               \langle a_s,a_{s+1} \rangle \end{array} \right), \quad
  \calE = (\calE_{jk}), \qquad
  j,k\in\{1,\ldots,s\}.
$$
The equation
$(I+\calE)^{-1} = I - \calE + \calE^2 - \calE^3 + \ldots$ and the inequality (\ref{aEe}) imply
$$
  \big|(I+\calE)^{-1}_{jk} - \delta_{jk}\big| \le 2\sigma_{jk}.
$$
Since $\lambda = (I+\calE)^{-1} \mu$, we have the estimate
$$
     |\lambda_j|
  =  \bigg| \sum_{k=1}^s (I+\calE)^{-1}_{jk} \mu_k
     \bigg|
 \le \sigma_{js+1} + 2\sum_{k=1}^s \sigma_{jk}\sigma_{ks+1}
 \le 2\sigma_{j\,s+1} .
$$
Putting $b=b'/|b'|$, we obtain by (\ref{|a-1|}):
\begin{eqnarray*}
&\displaystyle
       |b' - a_{s+1}|
  \le  \sum_{j=1}^s 2\sigma_{js+1} |a_j|
  \le  \sum_{j=1}^s 2\sigma_{js+1} (1 + \sigma_{jj}),  & \\
&\displaystyle
       \frac{|b'|}{|b|} = 1 + \sigma',  \quad
       |\sigma'|
   =   \big| |b'| - 1 \big|
  \le  \sigma_{s+1 s+1} + \sum_{j=1}^s  2\sigma_{j s+1} (1 + \sigma_{jj}).  &
\end{eqnarray*}
Lemma \ref{lem:GS} is proved.  \qed

We apply Lemma \ref{lem:GS} repeatedly, taking as $a_1,\ldots,a_{s+1}$ first,
$w_{-N},\ldots,w_{N-1},w_N$, then $w_{1-N},\ldots,w_{N-1},w_{-N}$, then $w_{1-N},\ldots,w_{N-1}$, and so on, and finally, $w_0$.

On each step we denote by $u_m$ the vector $b$ which replaces $w_m = a_{s+1}$. In this way we associate with the vectors $w_{-N},\ldots,w_N$ the orthonormal basis $u_{-N},\ldots,u_N$.

We have:
\begin{equation}
\label{um}
    u_m
  = (1 + \sigma')^{-1} \Big( w_m - \sum_{-|m|\le j < |m|,\, j\ne m} \lambda'_j w_j \Big).
\end{equation}
where the coefficients $\lambda'_{-|m|},\ldots,\lambda'_{|m|-1}$ are $\lambda_1,\ldots,\lambda_s$ with properly changed indices and
\begin{equation}
\label{|sigma|}
  |\sigma'| \le \eps_{mm} + \sum_{-|m|\le j < |m|,\, j\ne m} 2\eps_{jm} (1 + \eps_{jj}).
\end{equation}
For any $k\in \{-N,\ldots,N\}$ inequalities (\ref{um}), (\ref{|sigma|}), and (\ref{b-a}) imply the estimates
$$
      |u_{km} - w_{km}|
  \le \frac{\sigma'}{1 - \sigma'} |w_{km}|
     + \frac2{1-\sigma'} \sum_{-|m|\le j < |m|,\, j\ne m} \eps_{jm} |w_{kj}|.
$$
If $N$ is sufficiently large, we have: $\eps_{jj} < 1/2$ and $|\sigma'| < 1/2$. Thus
$$
      |u_{km} - w_{km}|
  \le 2\sigma' |w_{km}| + \sum_{-|m|\le j < |m|,\, j\ne m} 4\eps_{jm} |w_{kj}|.
$$
This proves (\ref{w'}). \qed

\subsection{Step 2}
\label{sec:M+M}

In this section we prove Theorem \ref{theo:mono1}. Take a big number $N$. Let
$U=(u_{jk})\in U(2n+1)$ be the operator constructed in Proposition \ref{prop:U}. Let
$\hat U = (\hat u_{jk}) \in U(2N+1)$ be an operator which diagonalizes $U$:
$$
  \hat U^{-1} U\hat U = D = \diag(D_{jj})_{-N\le j\le N}.
$$
Consider the unitary operator $\hat\calU$ on $l^2(\mZ)$ such that
$$
    \hat\calU_{jk}
  = \left\{ \begin{array}{cl}
              \hat u_{jk} & \mbox{ if both $|j|\le N$ and $|k|\le N$}, \\
              \delta_{jk} & \mbox{ if $|j| > N$ or $|k| > N$}.
            \end{array}
    \right.
$$
We define $\tilde\calN = \hat\calU^{-1} \calN\hat\calU = \hat\calU^* \calN\hat\calU$. It remains to check that $\tilde\calN - \tilde\calM^{(0)} - \tilde\calM^{(d)} = \tilde\calM^{(c)}$, where $\tilde\calM^{(0)}$ is determined by (\ref{tildeM}), $d_j=D_{jj}$ for $|j|\le N$,
$\tilde\calM^{(d)}$ is determined by (\ref{tildeMd}), and $\tilde\calM^{(c)}$ satisfies (\ref{Mc}).
\smallskip

We have: $\tilde N_{jk} = \Sigma_{uu} + \Sigma_{\delta u} + \Sigma_{u\delta} + \Sigma_{\delta\delta}$,
\begin{eqnarray*}
      \Sigma_{uu}
  &=& \sum_{-N\le j',k'\le N} \hat u^*_{jj'} \calN_{j'k'} \hat u_{k'k}, \\
      \Sigma_{\delta u}
  &=& \sum_{|j'|>N,\,-N\le k'\le N} \delta_{jj'} \calN_{j'k'} \hat u_{k'k}, \\
      \Sigma_{u\delta}
  &=& \sum_{-N\le j'\le N,\, |k'|>N} \hat u^*_{jj'} \calN_{j'k'} \delta_{k'k}, \\
      \Sigma_{\delta\delta}
  &=& \sum_{|j'|>N,\,|k'|>N} \delta_{jj'} \calN_{j'k'} \delta_{k'k}.
\end{eqnarray*}
There are 4 cases.
\smallskip

{\bf (a)}. If both $|j|>N$ and $|k|>N$, we have: $\Sigma_{uu}=\Sigma_{\delta u}=\Sigma_{u\delta}=0$. Hence $\tilde\calN_{jk} = \calN_{jk}$ and (\ref{tildeMc}) holds by (\ref{Mc}).
\smallskip

{\bf (b)}. If $|j|>N$ and $|k|\le N$ then $\Sigma_{uu}=\Sigma_{u\delta}=\Sigma_{\delta\delta}=0$, so that
$$
  \tilde\calN_{jk} = \sum_{-N\le k'\le N} \calN_{jk'} \hat U_{k'k}.
$$
We estimate this sum by using (\ref{Mc}) and the Cauchy inequality:
\begin{eqnarray*}
       |\tilde\calN_{jk}|^2
 &\le& \sum_{-N\le k'\le N} |\calN_{jk'}|^2  \sum_{-N\le k'\le N} |\hat U_{k'k}|^2 \\
 &\le& \sum_{-N\le k'\le N}
          \frac{c^2 e^{-2\beta |j-k'|}}
               {\langle j\rangle^2 \langle k'\rangle^2 \langle j-k'\rangle^{2\alpha-2}}
  \le   \frac{c^2 c_2 e^{-2\beta (|j|-N)}}{\langle j\rangle^4},
\end{eqnarray*}
where $c_2=c_\nu|_{\nu=2}$ is a constant from Lemma \ref{lem:tech1}.
\smallskip

{\bf (c)}. In the case $|j|\le N$ and $|k| > N$ the argument is the same as in {\bf (b)}.
\smallskip

{\bf (d)}. If both $|j|\le N$ and $|k|\le N$ then
$\Sigma_{\delta u}=\Sigma_{u\delta}=\Sigma_{\delta\delta}=0$. Then
\begin{eqnarray*}
      \tilde\calN_{jk}
  &=& \sum_{-N\le j',k'\le N} \hat U^*_{jj'} (u_{j'k'} + w'_{j'k'}) U_{k'k}
   =  D_{jj} \delta_{jk} + S, \\
      S
  &=& \sum_{-N\le j',k'\le N} \hat U^*_{jj'} w'_{j'k'} \hat U_{k'k}
\end{eqnarray*}
The coefficients $w'_{j'k'}$ satisfy (\ref{w'}). We put
$Q_{j'k}=\sum_{-N\le k'\le N} w'_{j'k'} \hat u_{k'k}$. By using the Cauchy inequality we have:
\begin{eqnarray*}
        | Q_{j'k} |^2
  &\le& \sum_{-N\le k'\le N} \frac{c_w^2}{\langle k'\rangle^2 (N+1)^4 \langle j'-k'\rangle^{2\alpha-2}}
   \le  \frac{c_w^2 c_2}{\langle j'\rangle^2 (N+1)^4}, \\
        |S|^2
  &=&  \bigg| \sum_{-N\le j'\le N} \hat u_{jj'} Q_{j'k} \bigg|^2
  \le  \sum_{-N\le j'\le N} \frac{c_w^2 c^2_2}{\langle j'\rangle^2 (N+1)^4}
  \le  \frac{c_w^2 c^2_2}{(N+1)^4} .
\end{eqnarray*}
This finishes proof of Theorem \ref{theo:mono1}.  \qed

\section{Technical statements}
\label{sec:tech}

\subsection{Several sums}

\begin{lem}
\label{lem:tech0}
1. For any $k,l\in\mZ$, where $k\ne\pm l$
\begin{equation}
\label{ineq1}
  \frac{1}{|k^2-l^2|} \le \frac{3}{\langle k\rangle + \langle l\rangle} .
\end{equation}
2. For integer $k_0,k,j$, where $j\ne\pm k_0$, $j\ne\pm k$
\begin{equation}
\label{ineq2}
  \frac{1}{|k_0^2-j^2| |k^2-j^2|} \le \frac{12}{\langle k_0\rangle^2} .
\end{equation}
\end{lem}

{\it Proof}. 1. W.l.g. $0\le l<k$. Then $k=l+s$, $s\ge 1$ and (\ref{ineq1}) is equivalent to the obvious inequality
$$
  2l+s+2 \le 6ls+3s^2, \qquad
  l\ge0, \quad
  s\ge1.
$$

2. W.l.g. $k_0,k$, and $j$ are nonnegative and $k_0\ge 1$. If $j\ge k_0/2$, then by (\ref{ineq1})
$$
      \frac1{|k_0^2-j^2| |k^2-j^2|}
  \le \frac{12}{(\langle k_0\rangle + \langle k_0/2\rangle) (\langle k\rangle + \langle k_0/2\rangle)}
  \le \frac{12}{\langle k_0\rangle^2}.
$$
If $0\le j\le k_0/2$, then
$$
      \frac1{|k_0^2-j^2|}
  \le \frac4{3 k_0^2}
  \le \frac{16}{3\langle k_0\rangle^2}.
$$
This implies (\ref{ineq2}). \qed

\begin{lem}
\label{lem:tech1}
Suppose $\nu>1$ and $\beta\ge 0$. Then
\begin{equation}
\label{tech1}
       \sum_{k=-\infty}^{+\infty}
         \frac{e^{-\beta|s-k|-\beta|k-m|}}
              {\langle s-k\rangle^\nu \langle k-m\rangle^\nu}
  \le  \frac{c_\nu e^{-\beta |s-m|}}{\langle s-m\rangle^\nu},
\end{equation}
where the constant $c_\nu$ depends only on $\nu$.
\end{lem}

{\it Proof}. Since $|s-k|+|k-m| \ge |s-m|$, it is sufficient to consider the case $\beta=0$.
Suppose $s\le m$. Then it is sufficient to split the sum into the following 4 pieces:
$$
  \sum_{k\le s}, \quad
  \sum_{s < k\le (s+m)/2}, \quad
  \sum_{(s+m)/2 < k\le m}, \quad
  \sum_{k > m}
$$
and to check estimate (\ref{tech1})$|_{\beta=0}$ for each of these sums.

The case $s>m$ is analogous. \qed

\subsection{Proof of Lemma \ref{lem:Gn}}
\label{sec:Gn}

Since $V\in C^{\alpha,\beta,1}$, we have:
$$
      |G_{\pm kk}|
   =  0, \quad   
      |G_{jk}|
  \le \frac{c_v e^{-\beta |j-k|}}{|j^2-k^2|\, \langle j-k\rangle^\alpha} \quad
      \mbox{for any integer $j\ne\pm k$}.
$$

Then by (\ref{ineq1})
$$
      |G_{jk}|
  \le \frac{3 c_v e^{-\beta |j-k|}}{(\langle j\rangle+\langle k\rangle)\, \langle j-k\rangle^\alpha}.
$$
By (\ref{tech1}) we have:
$$
       |(G^2)_{jk}|
   =   \Big| \sum_{s\in\mZ} G_{js} G_{sk} \Big|
  \le  \sum_{s\in\mZ}
        \frac{9 c^2_v e^{-\beta |j-k|}}
          {\langle j\rangle \langle k\rangle \langle j-s\rangle^\alpha \langle s-k\rangle^\alpha}
  \le   \frac{9 c^2_v c_\alpha e^{-\beta |j-k|}}
              {\langle j\rangle \langle k\rangle \langle j-k\rangle^\alpha }.
$$

Assuming that the second inequality (\ref{Gn<}) holds for some $n\ge 2$, we have:
$$
       |(G^{n+1})_{jk}|
   =   \Big| \sum_{s\in\mZ} (G^n)_{js} G_{sk} \Big|
  \le  \sum_{s\in\mZ}
       \frac{ (3c_v)^{n+1} c_\alpha^{n-1} e^{-\beta |j-k|}}
           {\langle j\rangle \langle k\rangle \langle j-s\rangle^\alpha \langle s-k\rangle^\alpha} \\
  \le  \frac{ (3c_v)^{n+1} c_\alpha^n e^{-\beta |j-k|}}
              {\langle j\rangle \langle k\rangle \langle j-k\rangle^\alpha }.
$$
\qed

\subsection{Proof of Lemma \ref{lem:Phi}}
\label{sec:Phi}

Below $c_0,c_1,\ldots$ are absolute constants.
\smallskip

$(\ba)$. Suppose $k=k_0$. Then by (\ref{psi1}) and (\ref{Phi<}) $\Phi_{k_0}(t) = A_1 + A_2 + A_3 + A_4$,
\begin{eqnarray*}
     A_1
 &=& \int_0^t \sum_{j\ne\pm k_0}
       \frac{v_{k_0-j}(\tau) v_{j-k_0}(\tau)}
            {k_0^2-j^2} \, d\tau , \\
     A_2
 &=& - \int_0^t \sum_{j\ne\pm k_0}
        \frac{v_{k_0-j}(\tau) v_{j-k_0}(0) e^{\frac{\tau}{ih} (j^2 - k_0^2)} }
            {k_0^2-j^2} \, d\tau , \\
\displaystyle
     A_3
 &=& - \int_0^t \sum_{j\ne\pm k_0} v_{k_0-j}(\tau) e^{\frac\tau{ih} (j^2-k_0^2)}
          \bigg( \int_0^\tau \frac{\dot v_{j-k_0}(\sigma)}{k_0^2-j^2}
                               e^{\frac\sigma{ih} (k_0^2-j^2)}\, d\sigma
          \bigg)\, d\tau, \\
     A_4
 &=& \int_0^t \frac{v_{2k_0}(\tau)}{ih}
                      \bigg( \int_0^\tau v_{-2k_0}(\sigma)\, d\sigma \bigg) d\tau .
\end{eqnarray*}

$(\ba_1)$. To estimate $A_1$, we have to use some cancelations. Putting $l= |j-k_0|$, we have:
\begin{eqnarray*}
      A_1
 &=&  \int_0^t \bigg(
                      \sum_{0 < l \ne 2|k_0|} v_{-l}(\tau) v_l(\tau)
                          \Big( \frac1{k_0^2 - (k_0+l)^2} + \frac1{k_0^2 - (k_0-l)^2} \Big)
                    + \frac{v_{-2k_0}(\tau) v_{2k_0}(\tau)}{k_0^2 - (3k_0)^2}
                             \bigg) \, d\tau \\
 &=&  \int_0^t \bigg(
                      \sum_{0 < l \ne 2|k_0|} v_{-l}(\tau) v_l(\tau)
                           \frac2{4k_0^2 - l^2}
                    - \frac{v_{-2k_0}(\tau) v_{2k_0}(\tau)}{8k_0^2}
                             \bigg) \, d\tau
\end{eqnarray*}
This implies
\begin{equation}
\label{A}
       |A_1|
  \le  |t| \bigg( \sum_{0<l\ne 2|k_0|}
                     \frac{2 c_v^2 e^{-2\beta |l|}}{\langle l\rangle^{2\alpha} |4k_0^2 - l^2|}
                   + \frac{c_v^2 e^{-4\beta |k_0|}}{\langle 2k_0\rangle^{2\alpha} 8k_0^2}
         \bigg)
  \le  \frac{|t| c_v^2 c_0}{\langle k_0\rangle^2} .
\end{equation}

$(\ba_2)$. We estimate $A_2$ by using integration by parts and (in the last inequality) by estimate (\ref{ineq1}):
\begin{eqnarray}
\nonumber
      |A_2|
  &=& \bigg| \int_0^t \sum_{j\ne\pm k_0} \frac{h v_{k_0-j}(\tau) v_{j-k_0}(0)}{(k_0^2-j^2)^2}\,
                   d e^{\frac{\tau}{ih}(j^2 - k_0^2)} \bigg|
\,\le\, \sum_{j\ne\pm k_0} \frac{(2+|t|) h c_v^2 e^{-2\beta |j-k_0|}}
                                {\langle k_0-j\rangle^{2\alpha}(k_0^2-j^2)^2} \\
\label{B}
 &\le& \frac{(2+|t|) h c_v^2 c_1}{\langle k_0\rangle^2}.
\end{eqnarray}

$(\ba_3)$. The equation
$$
       A_3
  =    \int_0^t \sum_{j\ne\pm k_0}
          \frac{ih v_{k_0-j}(\tau)}{k_0^2 - j^2}
          \bigg( \int_0^\tau \frac{\dot v_{j-k_0}(\sigma)}{k_0^2-j^2}
                               e^{\frac\sigma{ih} (k_0^2-j^2)}\, d\sigma
          \bigg)\, d e^{\frac\tau{ih} (j^2-k_0^2)},
$$
implies that
\begin{equation}
\label{B'}
       |A_3|
  \le  \sum_{j\ne\pm k_0} \frac{(3|t| + t^2/2) h c_v^2}
                               {(k_0^2-j^2)^2 \langle k_0-j\rangle^{2\alpha}}
  \le  \frac{(3|t| + t^2/2) h c_2 c_v^2}{\langle k_0\rangle^2}.
\end{equation}

$(\ba_4)$. The estimate
\begin{equation}
\label{C}
   |A_4| \le \frac{t^2 c_v^2 c_2}{h \langle k_0\rangle^2}
\end{equation}
is obvious.
\smallskip

$(\bb)$. Now suppose $k\ne\pm k_0$. Then $\Phi_k(t) = B_1+B_2+B_3+B_4$,
\begin{eqnarray*}
     B_1
 &=& \int_0^t \sum_{j\ne\pm k_0}
        \frac{v_{k-j}(\tau) v_{j-k_0}(\tau) e^{\frac\tau{ih} (k_0^2-k^2)}}{k_0^2 - j^2} \, d\tau, \\
     B_2
 &=& - \int_0^t \sum_{j\ne\pm k_0}
        \frac{v_{k-j}(\tau) v_{j-k_0}(0) e^{\frac\tau{ih} (j^2-k^2)}}{k_0^2 - j^2} \, d\tau, \\
     B_3
 &=& - \int_0^t \sum_{j\ne\pm k_0} v_{k-j}(\tau) e^{\frac\tau{ih} (j^2-k^2)}
          \bigg( \int_0^\tau \frac{\dot v_{j-k_0}(\sigma)}{k_0^2-j^2}
                               e^{\frac\sigma{ih} (k_0^2-j^2)}\, d\sigma
          \bigg)\, d\tau, \\
     B_4
 &=& \int_0^t
        \frac{v_{k+k_0}(\tau)}{ih} e^{\frac\tau{ih} (k_0^2-k^2)}
        \Big( \int_0^\tau v_{-2k_0}(\sigma)\, d\sigma \Big) \, d\tau.
\end{eqnarray*}

$(\bb_1)$. Since $k^2\ne k_0^2$, we can use in $B_1$ integration by parts:
$$
    B_1
  = \int_0^t \sum_{j\ne\pm k_0}
        \frac{ih v_{k-j}(\tau) v_{j-k_0}(\tau)}{(k_0^2-k^2)(k_0^2 - j^2)}
            \, d e^{\frac\tau{ih} (k_0^2-k^2)} .
$$
Then Lemmas \ref{lem:tech1} and \ref{lem:tech0} imply the estimate
\begin{equation}
\label{D}
      |B_1|
 \le  \sum_{j\ne\pm k_0}
        \frac{(2 + |t|) h c_v^2 e^{-\beta |k-k_0|}}
             {(k_0^2-k^2)(k_0^2-j^2) \langle k-j\rangle^\alpha \langle j-k_0\rangle^\alpha} 
 \le    \frac{9(2 + |t|) h c_v^2 c_\alpha  e^{-\beta |k-k_0|}}
             {\langle k_0\rangle^2 \langle k-k_0\rangle^\alpha}.
\end{equation}

$(\bb_2)$. In the sum $B_2$ the term with $j=k$ vanishes because $v_0=0$. However we have to consider separately the term with $j=-k$. We denote this term $B'_2$. Integrating by parts, we obtain:
$$
     B_2 - B'_2
  =  \int_0^t \sum_{j\ne\pm k_0,\, j\ne\pm k}
       \frac{ih v_{k-j}(\tau) v_{j-k_0}(0)}{(k_0^2 - j^2)(k^2 - j^2)} \, d e^{\frac\tau{ih} (j^2-k^2)}.
$$
This implies the estimate
$$
     |B_2 - B'_2|
 \le \sum_{j\ne\pm k_0,\, j\ne\pm k}
        \frac{(2 + |t|) h c_v^2 e^{-\beta |k-k_0|}}
             {|k_0^2 - j^2| |k^2 - j^2| \langle k-j\rangle^\alpha \langle j-k_0\rangle^\alpha} .
$$
By Lemmas \ref{lem:tech0} and \ref{lem:tech1}
\begin{equation}
\label{E}
     |B_2 - B'_2|
 \le \sum_{j\ne\pm k_0,\, j\ne\pm k}
        \frac{12 (2 + |t|) h c_v^2 c_\alpha e^{-\beta |k-k_0|}}
             {\langle k_0\rangle^2 \langle k-k_0\rangle^\alpha} .
\end{equation}
We have:
\begin{equation}
\label{F}
     |B'_2|
  =   \bigg| \int_0^t \frac{v_{2k}(\tau) v_{-k-k_0}(0)}{k_0^2 - k^2} \, d\tau
      \bigg|
 \le  \frac{|t| c_v^2 e^{-\beta |k-k_0|}}
              {|k_0^2 - k^2| \langle 2k\rangle^\alpha \langle k+k_0\rangle^\alpha}
 \le  \frac{|t| c_v^2 c_3 e^{-\beta |k-k_0|}}
            {\langle k_0\rangle^2 \langle k-k_0\rangle^{\alpha-1}}
\end{equation}
(the worst situation is $k=-k_0\pm 1$).
\smallskip

$(\bb_3)$. We have:
$$
       B_3
  =    \int_0^t \sum_{j\ne\pm k_0} \frac{ih v_{k-j}(\tau)}{k^2-j^2}
          \bigg( \int_0^\tau \frac{\dot v_{j-k_0}(\sigma)}{k_0^2-j^2}
                               e^{\frac\sigma{ih} (k_0^2-j^2)}\, d\sigma
          \bigg)\, d e^{\frac\tau{ih} (j^2-k^2)}
$$
Therefore by Lemmas \ref{lem:tech0} and \ref{lem:tech1}
\begin{equation}
\label{F'}
      |B_3|
  \le \sum_{j\ne\pm k_0}
        \frac{(3|t|+t^2/2) h c_v^2 e^{-\beta |k-k_0|}}
             {|k^2 - j^2| |k_0^2 - j^2| \langle k-j\rangle^\alpha \langle j-k_0\rangle^\alpha}
  \le   \frac{(3|t|+t^2/2) h c_v^2 c_4 e^{-\beta |k-k_0|}}
             { \langle k_0\rangle^2 \langle k-k_0\rangle^\alpha}
\end{equation}

$(\bb_4)$. Now consider $B_4$:
$$
     B_4
  =  \int_0^t  \frac{v_{k+k_0}(\tau)}{k_0^2 - k^2}
               \bigg( \int_0^\tau v_{-2k_0}(\sigma)\, d\sigma \bigg)\, d e^{\frac\tau{ih}(k_0^2-k^2)} .
$$
This implies
\begin{equation}
\label{G}
     |B_4|
 \le \frac{(3|t|+t^2/2) c_v^2 e^{-\beta |k-k_0|}}
          { |k_0^2 - k^2| \langle k+k_0\rangle^\alpha \langle 2k_0\rangle^\alpha}
 \le \frac{(3|t|+t^2/2) c_v^2 c_5  e^{-\beta |k-k_0|}}
          {\langle k_0\rangle^2 \langle k-k_0\rangle^{\alpha-1}}.
\end{equation}
(the worst situation is again $k=-k_0\pm 1$).

$(\bc)$. Finally consider the case $k=-k_0$. Then $\Phi_{-k_0}(t) = C_1+C_2$, where
\begin{eqnarray*}
      C_1
 &=&  \int_0^t \sum_{j\ne \pm k_0}
       \frac{v_{-k_0-j}(\tau) v_{j-k_0}(\tau)
            - v_{-k_0-j}(\tau) v_{j-k_0}(0) e^{\frac\tau{ih}(j^2-k_0^2)}}
            {k_0^2 - j^2} \, d\tau, \\
      C_2
 &=&  \int_0^t \sum_{j\ne\pm k_0} v_{-k_0-j}(\tau) e^{\frac{\tau}{ih}(j^2-k_0^2)}
               \bigg( \int_0^\tau \frac{\dot v_{j-k_0}(\sigma)}{k_0^2-j^2} 
                                 e^{\frac\sigma{ih} (k_0^2-j^2)}\, d\sigma 
               \bigg) \, d\tau .
\end{eqnarray*}

$(\bc_1)$. We have:
$$
       |C_1|
  \le  \sum_{j\ne\pm k_0}
         \frac{2c_v^2 |t|  e^{-2\beta |k_0|}}
              {|k_0^2-j^2| \langle k_0+j\rangle^\alpha \langle k_0-j\rangle^\alpha}
  \le  \frac{4c_v^2 |t|}{\langle k_0\rangle} \sum_{j\ne\pm k_0}
         \frac{e^{-2\beta |k_0|}}
              {\langle k_0+j\rangle^\alpha \langle k_0-j\rangle^\alpha} .
$$
By Lemma \ref{lem:tech1}
\begin{equation}
\label{I}
      |C_1|
  \le \frac{4|t| c_v^2 c_\alpha e^{-2\beta |k_0|}}
           {\langle k_0\rangle \langle 2k_0\rangle^\alpha}
  \le \frac{4|t| c_v^2 c_\alpha  e^{-2\beta |k_0|}}
           {\langle k_0\rangle^2 \langle 2k_0\rangle^{\alpha-1}}.
\end{equation}

$(\bc_2)$. We may estimate $C_2$ without integration by parts:
\begin{equation}
\label{J}
      |C_2|
  \le \sum_{j\ne\pm k_0}
          \frac{t^2 c_v^2 e^{-2\beta |k_0|}}
               {2 |k_0^2-j^2| \langle k_0+j\rangle^\alpha \langle k_0-j\rangle^\alpha}
  \le     \frac{t^2 c_v^2 c_6 e^{-2\beta |k_0|}}
               { \langle k_0\rangle^2 \langle k_0\rangle^{\alpha-1} }.
\end{equation}

Lemma \ref{lem:Phi} follows from the estimates (\ref{A})--(\ref{J}).  \qed

\end{document}